\documentclass[11pt]{article} 
\usepackage[utf8]{inputenc}
\usepackage{caption,graphicx,amsmath,amsthm,amssymb, authblk,xcolor,enumitem,subfigure,float,verbatim,hyperref,pgfplots,bibentry,mathtools,soul,isomath, cancel}
\usepackage{mathrsfs}
\usepackage{indentfirst}
\usepackage{nccmath}
\graphicspath{ {images/} }
\usepackage{multicol, MnSymbol}
\usepackage{graphicx,amsmath,amssymb,amsthm,marvosym,tikz,fontenc}
\usetikzlibrary{backgrounds}
\usetikzlibrary{patterns, shapes.misc, positioning}
\usepackage{pgfplots}
\pgfplotsset{width=10cm,compat=1.9,tick scale binop=\times}
\usepackage{indentfirst}
\usepackage[a4paper,bindingoffset=0.2in,%
left=0.8in,right=1in,top=1.2in,bottom=1.2in,%
footskip=.25in]{geometry}
\usepackage{relsize}
\usepackage[english]{babel}
\allowdisplaybreaks
\theoremstyle{plain}
\newtheorem{theorem}{Theorem}[section]
\newtheorem{lemma}[theorem]{Lemma}
\newtheorem{proposition}[theorem]{Proposition}

\newtheorem{definition}[theorem]{Definition}

\newtheorem{example}[theorem]{Example}

\date{}
\begin{document}
\title
{\bf{On the Construction of New Families of Distance Equienergetic Graphs}}
\author {\small Anjitha Ashokan
\footnote{anjithaashokan1996@gmail.com}   and Chithra A V
\footnote{chithra@nitc.ac.in} \\ \small Department of Mathematics, National Institute of Technology Calicut,\\
\small Kerala, India-673601}
\date{}
\maketitle
\begin{abstract}
    For a graph $G$ with vertex set $V(G)=\{v_{i}:1\leq i \leq p\}$ and edge set $E(G)=\{e_{j}:1\leq j \leq q\},$  its subdivision graph $S_G$ is generalized to $S_{G}(n,m),$ by replacing each vertex $v_{i}$
    with $n$ vertices and each edge-vertex $e_{j}$ with $m$ vertices,
with each copy having the same neighbors as the vertex it replaces.
 Building on this construction, we define the generalized subdivision vertex join $S_{G_{1}}(n,m)\dot{\vee} G_{2}$ and the generalized subdivision edge join  $S_{G}(n,m)\veebar G_{2},$ which extend the known subdivision vertex and subdivision edge join of graphs. In this paper, we determine the distance spectrum of $S_{G_{1}}(n,n)\dot{\vee} G_{2}$ and $S_{G}(n,n)\veebar G_{2},$ when $G_{1}$ and $G_{2}$ are regular graphs.
    As an application, we construct new families of $D$-cospectral graphs and $D$-equienergetic graphs of diameter $4.$
    
\end{abstract}
\section{Introduction}

Let $G$ be a simple connected graph with vertex set $V(G)=\{v_{i}:1\leq i \leq p\}$ and edge set $E(G)=\{e_{j}:1\leq j \leq q\}.$ 
The \textit{adjacency matrix} of  $G,$  $A(G)=(a_{ij})$ is a $p\times p$  symmetric matrix such that $a_{ij}=1,$ if $v_{i}$ and $v_{j}$ are adjacent, and $0$ otherwise.
Let  $\lambda_{1}, \lambda_{2},\ldots, \lambda_{p}$ be the $eigenvalues$ of $A(G)$ (adjacency eigenvalues of $G$). 
The \textit{adjacency spectrum} of $G$, $spec (G)$, is the collection of all eigenvalues of $A(G)$ along with their multiplicities.
Spectral graph theory provides a powerful framework for investigating the structural and  combinatorial properties of graphs through the eigenvalues of matrices associated with them. Among the most commonly studied matrices are the adjacency matrix, the Laplacian matrix, and the distance matrix, each of which captures different aspects of a graph's structure. In particular, the distance matrix has attracted considerable attention owing to its applications in communication networks, chemistry, and various optimization problems; see \cite{zbMATH06319411}.

The \textit{distance} between  any two vertices $v_{i},v_{j}\in V(G),$  $d(v_{i},v_{j})$ is defined as the length of the shortest path connecting them in $G.$
In 1971, Graham and Pollack \cite{zbMATH03363442}  introduced the distance matrix of a connected graph to study data communication problems.
The \textit{distance matrix} of $G,$ $D(G)$ is the $p\times p$ symmetric matrix whose $(i, j)$-th entry is equal to $d(v_{i},v_{j}).$
The eigenvalues of $D(G)$ are called the $D$-\textit{eigenvalues} of $G$ and they along with their multiplicities form the $D$-\textit{spectrum} of $G.$
The {$D$-eigenvalues of $G$ }have applications in various fields like data communication, chemistry, and molecular biology (see \cite{balasubramanian1990computer} and the references therein).
Let $\eta_{i}$ be  the $i^{th}$ largest eigenvalue of $D(G)$.
The concept of distance energy of a graph was introduced in 2008 \cite{indulal2008distance}.
The \textit{distance energy} ($D$-energy) of a graph $G$ is defined as the
sum of the absolute values of the eigenvalues of its
distance matrix.

 The subdivision graph  $S_G$ of a graph $G$ \cite{cvetkovic1980spectra}, is obtained from $G$ by subdividing each edge of $G$ using a new vertex,
 the vertex set of $S_{G}$ is $V(G)\cup E(G).$
The concept of subdivision vertex join (respectively, subdivision edge join), which is obtained from $S_{G_{1}}$ and $G_{2}$ by joining each vertex of $V(G_{1})$ (respectively, $E(G_{1})$) to every vertex of $V(G_{2}),$ was introduced in 2012 \cite{zbMATH06221129}.
Later in 2024, Shamsher and Pirzada \cite{zbMATH08041341} generalized the notion of subdivision graphs as follows:

 Let $G$ be a graph with vertex set $V(G)=\{v_{1},v_{2},\ldots,v_{p}\}$ and edge set $E(G)=\{e_{1},e_{2},\ldots,e_{q}\}.$
 The generalized subdivision graph of $G$, denoted by
$S_G(n,m)$, is the graph obtained from $S_G$ by replacing
each vertex $v_i$, $i=1,2,\ldots,p$, with $n$ new vertices
having the same neighborhood as $v_i$, followed by replacing
each edge $e_j$, $j=1,2,\ldots,q$, with $m$ new vertices
adjacent to the neighbors of $e_j$.
Equivalently, the graph $S_{G}(n,m)$ has vertex set $$\Big(V(G)\times \{1,2,\ldots,n\}\Big)\cup \Big(E(G)\times \{1,2,\ldots ,m\}\Big),$$
where two vertices $(v_{l},i)$ and $(e_{m},j),$  are adjacent whenever the edge $e_{m}\in E(G)$ is incident with the vertex $v_{l}\in V(G)$ in $G.$ 
Motivated by this generalized subdivision graph, we introduce the concept of generalized subdivision vertex join and generalized subdivision edge join of two graphs $G_{1}$ and $G_{2}$ as follows.
\begin{definition}
    The generalized subdivision vertex join $S_{G_{1}}(n,m) \dot{\vee} G_{2}$ is the graph obtained from $S_{G_{1}}(n,m)$ and $G_{2}$ by joining each vertex of $V(G_{1})\times \{1,2,\ldots,n\}$ to every vertex of $V(G_{2})$. 
\end{definition}
Note that, if $G_{i}$ is a graph on $p_{i}$ vertices and $q_{i}$ edges, $i=1,2,$  then  $S_{G_{1}}(n,m) \dot{\vee} G_{2}$  has $p_{1}n+q_{1}m+p_{2}$ vertices and $p_{2}+p_{1}np_{2}+(p_{1}+q_{1})mn$ edges.
\begin{example}
   Let $G_{1}=K_{2}$ and $G_{2}=C_{3}.$ Then, $S_{K_{2}}(2,2) \dot{\vee} C_{3}$ is depicted in Figure \ref{fig:1}.
\end{example}
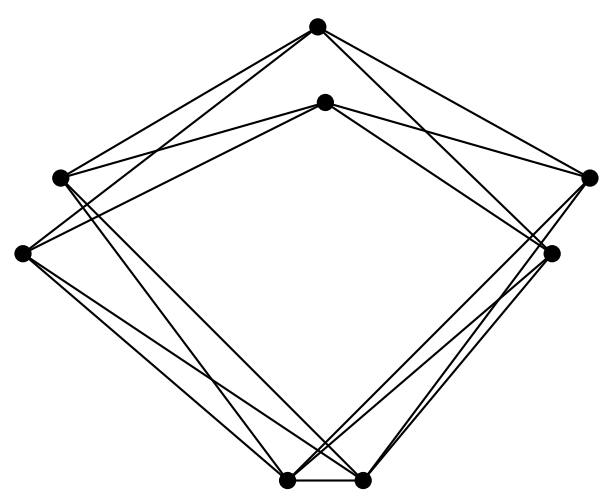
\begin{figure}[H]
 \centering
 \begin{tikzpicture}[scale=.5]
 \filldraw[fill=black](-20,1)circle(0.1cm); 
 \filldraw[fill=black](-20.5,0)circle(0.1cm); 
 \filldraw[fill=black](-13,1)circle(0.1cm); 
 \filldraw[fill=black](-13.5,0)circle(0.1cm); 
 \filldraw[fill=black](-16.5,2)circle(0.1cm); 
 \filldraw[fill=black](-16.6,3)circle(0.1cm);
 \filldraw[fill=black](-17,-3)circle(0.1cm); 
 \filldraw[fill=black](-16,-3)circle(0.1cm);
 \draw(-20, 1)--(-16.5,2);
 \draw(-20, 1)--(-16.6,3);
 \draw(-20, 1)--(-17,-3);
 \draw(-20, 1)--(-16,-3);
 \draw(-20.5, 0)--(-16.5,2);
 \draw(-20.5, 0)--(-16.6,3);
 \draw(-20.5, 0)--(-17,-3);
 \draw(-20.5, 0)--(-16,-3);
 \draw(-13, 1)--(-16.5,2);
 \draw(-13, 1)--(-16.6,3);
 \draw(-13, 1)--(-17,-3);
 \draw(-13, 1)--(-16,-3);
 \draw(-13.5, 0)--(-16.5,2);
 \draw(-13.5, 0)--(-16.6,3);
 \draw(-13.5, 0)--(-17,-3);
 \draw(-13.5, 0)--(-16,-3);
   \draw(-17,-3)--(-16,-3);
 \end{tikzpicture}
  \caption{$S_{K_{2}}(2,2)\dot{\vee}K_{2}$}
 \label{fig:1}
 \end{figure}

\begin{definition}
     The generalized subdivision edge join $S_{G_{1}}(n,m) \veebar G_{2}$ is the graph obtained from $S_{G_{1}}(n,m)$ and $G_{2}$ by joining each vertex of $E(G_{1})\times \{1,2,\ldots,m\}$ to every vertex of $V(G_{2})$. 
\end{definition}
Note that, if $G_{i}$ is a graph on $p_{i}$ vertices and $q_{i}$ edges, $i=1,2,$  then  $S_{G_{1}}(n,m) \veebar G_{2}$  has $p_{1}n+q_{1}m+p_{2}$ vertices and $p_{2}+q_{1}mp_{2}+(p_{1}+q_{1})mn$ edges.
\begin{example}
   Let $G_{1}=K_{2}$ and $G_{2}=C_{3}.$ Then, $S_{K_{2}}(2,2) \veebar C_{3}$ is depicted in Figure \ref{fig:2}.
\end{example}
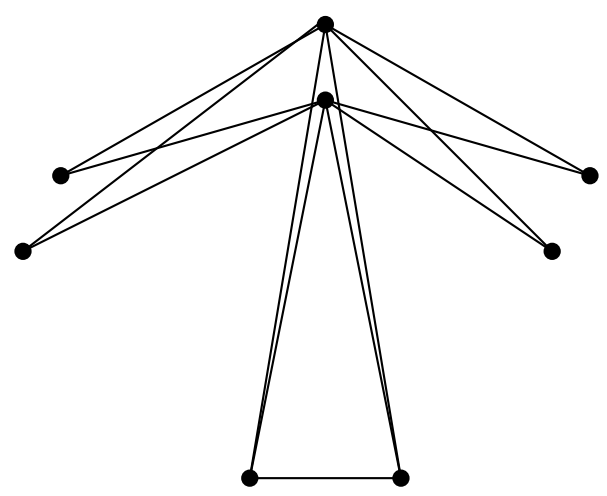
\begin{figure}[H]
 \centering
 \begin{tikzpicture}[scale=.5]
 \filldraw[fill=black](-20,1)circle(0.1cm); 
 \filldraw[fill=black](-20.5,0)circle(0.1cm); 
 \filldraw[fill=black](-13,1)circle(0.1cm); 
 \filldraw[fill=black](-13.5,0)circle(0.1cm); 
 \filldraw[fill=black](-16.5,2)circle(0.1cm); 
 \filldraw[fill=black](-16.5,3)circle(0.1cm);
 \filldraw[fill=black](-17.5,-3)circle(0.1cm); 
 \filldraw[fill=black](-15.5,-3)circle(0.1cm);
 \draw(-20, 1)--(-16.5,2);
 \draw(-20, 1)--(-16.5,3);
 \draw(-20.5, 0)--(-16.5,2);
 \draw(-20.5, 0)--(-16.6,3);
 \draw(-13, 1)--(-16.5,2);
 \draw(-13, 1)--(-16.5,3);
 \draw(-13.5, 0)--(-16.5,2);
 \draw(-13.5, 0)--(-16.5,3);
    \draw(-16.5,2)--(-17.5,-3);
  \draw(-16.5,2)--(-15.5,-3);  
      \draw(-16.5,3)--(-17.5,-3);
  \draw(-16.5,3)--(-15.5,-3);  
    \draw(-17.5,-3)--(-15.5,-3);
 \end{tikzpicture}
  \caption{$S_{K_{2}}(2,2)\veebar K_{2}$}
 \label{fig:2}
 \end{figure}

Two non-isomorphic graphs of the same order are cospectral  (respectively, $D$- cospectral) if they have the same adjacency spectrum (respectively, $D$- spectrum).
Non-isomorphic, non-cospectral (respectively, non $D$-cospectral) graphs of the same order are equienergetic (respectively, distance equienergetic) if their adjacency energies (respectively, distance energies) are the same.
The construction of non-isomorphic equienergetic graphs with respect to various graph matrices is an important problem in spectral graph theory. In particular, several constructions of equienergetic and distance equienergetic graphs have been investigated in the literature; see \cite{zbMATH08090417, zbMATH07833729, zbMATH05664447, zbMATH07833946, zbMATH07249081}. 
In this work, we investigate the distance spectra of generalized subdivision vertex (respectively, edge) join  of regular graphs. In particular, we derive new spectral results that encompass the main results of Indulal et al. (2019) \cite{zbMATH07090170} as special cases.
In this work, 
we study the distance spectra of generalized subdivision vertex (respectively, edge) join of regular graphs.
Our results generalize the main results of Indulal et al. (2019) \cite{zbMATH07090170}.
 Furthermore, using these results, we construct families of non-$D$-cospectral $D$-equienergetic graphs with diameter $4.$

This article is organized as follows: Section 2 presents the preliminaries and notations used throughout.
In Section 3, we present the distance spectra of the class of graphs obtained by the generalized subdivision vertex (respectively, edge) join of regular graphs.   Finally, in Section 4, we  construct  families of $D$-equienergetic and $D$-cospectral graphs of diameter 4.
\section{Preliminaries}
The $incidence$ $matrix$ $R=R(G)=(r_{ij})$ of  $G$, is the $p\times q$ matrix  such that $r_{ij}=1$ if and only if vertex $v_{i}$ and edge $e_{j}$ are incident in $G.$
The $line$ $graph$ $L(G)$ of $G$ is the graph with its vertex set as the edge set of $G$, and two vertices in $L(G)$ become adjacent if the corresponding edges of $G$ share a common end vertex. The adjacency matrix of $L(G)$ is denoted by $B(G).$
 We indicate by $K_{p}$ the complete graph, by $C_{p}$ the cycle, each of order $p.$
Throughout this paper, $J_{p\times q}$ and $I_{p\times p}$
represent the all-ones matrix and the identity matrix of the
indicated orders, respectively.
The Kronecker product $A\otimes B,$ of two matrices $A=(a_{ij})_{n\times m}$ and $B=(b_{ij})_{p\times q}$ is obtained by replacing the $(i,j)$-th entry in $A$ by $a_{ij}B.$ For matrices $A,B,C,$ and $D$ of appropriate orders, we have $(A\otimes B)(C\otimes D)=A C \otimes B D.$ Also, for convenience, we denote $I_{p\times p}\otimes (J-I)_{n \times n}$ as $\widetilde{I}_{pn\times pn}.$ The  \textit{equitable quotient matrix} \cite{brouwer2011spectra} of a partitioned matrix is denoted by M.\\
The following lemmas provide the foundation for the results developed in the subsequent sections.
\begin{lemma}$\textnormal{\cite{cvetkovic1980spectra}}$\label{thelinegraph}
    Let $G$ be an $r$\textcolor{red}{-}regular, $(p,q)$ graph with an adjacency matrix $A,$ and the incidence matrix $R$. Let $L(G)$ be the line graph of $G$. Then, $RR^{T}=A+rI$, $R^{T}R=B+2I$, where $B$ is the adjacency matrix of $L(G)$, also $JR=2J=R^{T}J$ and $JR^{T}=rJ=RJ.$ 
\end{lemma}
\begin{lemma}$\textnormal{\cite{cvetkovic1980spectra}}$\label{Beigenvalue}
    Let $G$ be an $r$-regular, $(p,q)$ graph with $spec(G)=\{r,\lambda_{2},\ldots, \lambda_{p}\}$. Then, 
                   \begin{equation*}
                       spec(L(G))=\begin{pmatrix}
                           2r-2 & \lambda_{2}+r-2 & \cdots&\lambda_{p}+r-2& -2\\
                           1 & 1 &\cdots& 1 & q-p
                       \end{pmatrix}.
                   \end{equation*}
 Also, $Z$ is an eigenvector corresponding to the eigenvalue $-2$ if and only if $RZ=0$.
 \end{lemma}
 
\section{Distance Spectra of Generalized Subdivision Vertex (Edge) Joins}
In this section, we determine the distance spectra of the generalized subdivision vertex, and edge join of regular graphs $G_{1}$ and $G_{2}$ in terms of their adjacency spectra and the adjacency spectrum of the line graph of $G_{1}.$

The following lemma establishes the distance between  vertices
in $S_{G_{1}}(n,n) \dot{\vee} G_{2}.$
\begin{proposition}
    Let $V(G)=\{v_{1},v_{2},\ldots,v_{p}\}$ and $E(G)=\{e_{1},e_{2},\ldots, e_{q}\}.$ Then, the distance between any two distinct vertices, say $x$ and $y$, in $S_{G_{1}}(n,n) \dot{\vee} G_{2}$ is as follows:
\begin{equation*}
   d(x,y)=\begin{cases}
            2 & \text{ if } x=(v_{i},j),y=(v_{l},s), \\
            1 &\text{ if } x= (v_{i},j), y=(e_{l},s) \text{ and } v_{i} \text{ is incident with } e_{l} \text{ in } G_{1}, \\
             3 &\text{ if } x= (v_{i},j), y=(e_{l},s) \text{ and } v_{i} \text{ is not incident with } e_{l} \text{ in } G_{1}, \\
             1& \text{ if } x=(v_{i},j)  \text{ and } y\in V(G_{2}),\\
             2& \text{ if } x= (e_{i},j), y=(e_{l},s) \text{ and } e_{i}, e_{l} \text{ have a common end point in } G_{1},\\
              4& \text{ if } x= (e_{i},j), y=(e_{l},s) \text{ and } e_{i}, e_{l} \text{ have  no common end point in } G_{1},\\
              2 &\text{ if } x=(e_{i},j)  \text{ and } y\in V(G_{2}),\\
              1  &\text{ if }x,y\in V(G_{2}) \text{ and }x,y \text{ are adjacent in }V(G_{2}),\\
              2  &\text{ if }x,y\in V(G_{2}) \text{ and }x,y \text{ are not adjacent in }V(G_{2}). 
            \end{cases}
\end{equation*}

\end{proposition}
\begin{lemma}\label{distance spectrum of generalized subdivision vertex join}
 Let $G_{i}$ be an $r_{i}$-regular graph on $p_{i}$ vertices and $q_{i}$ edges with the adjacency matrix $A(G_{i})$ and adjacency spectrum $\{r_{i},\lambda_{i2},\lambda_{i3},\ldots,\lambda_{ip_{i}}\},\text{ } i=1,2.$ Then, the distance spectrum of $S_{G_{1}}(n,n) \dot{\vee} G_{2}$ consists of
  $-2$ with multiplicity $(n-1)(p_{1}+q_{1});$
 $2(n-1)$ with multiplicity $q_{1}-1;$
 $-2(n\lambda_{1j}+nr_{1}+1), j=2,3,\ldots,p_{1}$;
     $-2-\lambda_{2m},m=2,3,\ldots,p_{2}$;
     and the eigenvalues of the matrix, $$\begin{pmatrix}
         2(np_{1}-1)&n(3q_{1}-2r_{1})&p_{2}\\
         n(3p_{1}-4)&4(q_{1}-r_{1})n+2(n-1)&2p_{2}\\
         p_{1}n&2q_{1}n&2(p_{2}-1)-r_{2}
        \end{pmatrix}.$$
\end{lemma}
\begin{proof}
   By suitable labeling of the vertices of $S_{G_{1}}(n,n) \dot{\vee} G_{2}$, the distance matrix $D(S_{G_{1}}(n,n) \dot{\vee} G_{2})$ can be represented as
   $$\scriptsize{\begin{pmatrix}
       \big(2(J-I)_{p_{1}\times p_{1}}\otimes J_{n\times n} \big)+2 \widetilde{I}_{p_{1}n\times p_{1}n}   & (3J-2R(G_{1}))_{p_{1}\times q_{1}} \otimes J_{n\times n}   &  J_{p_{1}n\times p_{2}}\\
       \big( (3J-2R(G_{1})^{T})_{q_{1}\times p_{1}}\otimes J_{n\times n}  \big)  & \big( (4(J-I)-2B(G_{1}))_{q_{1}\times q_{1}}\otimes J_{n\times n}\big)+2\widetilde{I}_{q_{1}n\times q_{1}n} & 2J_{q_{1}n\times p_{2}}\\
       J_{p_{2}\times p_{1}n} &2J_{p_{2}\times q_{1}n} &\big( 2(J-I)-A(G_{2})\big)_{p_{2}\times p_{2}}
   \end{pmatrix}}.$$
   
 For $s=1,2,3,\ldots, p_{1},$ consider the vector 
\begin{equation}\label{eq:eigenvector_gsvj one}
\begin{split}
 z^{(s,t)}=&\bigl(x_{1,1},x_{1,2},x_{1,3},\ldots, x_{1,n},
         x_{2,1},x_{2,2},x_{2,3},\ldots, x_{2,n},
         \cdots,
         x_{p_1,1},x_{p_1,2},x_{p_1,3},\ldots, x_{p_1,n},\\
        & y_{1,1},y_{1,2},y_{1,3},\ldots, y_{1,n},
          y_{2,1},y_{2,2},y_{2,3},\ldots, y_{2,n},
          \cdots,
          y_{q_1,1},y_{q_1,2},y_{q_1,3},\ldots, y_{q_1,n},\\
        & w_{1,1},w_{1,2},w_{1,3},\ldots,w_{1,p_{2}}\bigr) \in \mathbb{R}^{N},
\end{split}
\end{equation}
 where $N=p_{1}n+q_{1}n+p_{2},$ and 
 $$x_{s,1}=1 \text { and } x_{s,t}=-1, \text { for } t=2,3,\ldots, n$$ 
 and all remaining coordinates of $z^{(s,t)}$ are equal to $0$.
 Then, $$D(S_{G_{1}}(n,n) \dot{\vee} G_{2})z^{(s,t)}=-2z^{(s,t)}, s=1,2,3,\ldots,p_{1}. $$ 
This construction yields $(n-1)p_{1}$ linearly independent eigenvectors of $D(S_{G_{1}}(n,n) \dot{\vee} G_{2})$  associated with the eigenvalue $-2.$

 Using  similar argument,  for $k=1,2,3,\ldots, q_{1},$ consider the vector 
\begin{equation}\label{eq:eigenvector_gsvj two}
\begin{split}
 h^{(k,l)}=&\bigl(u_{1,1},u_{1,2},u_{1,3},\ldots, u_{1,n},
         u_{2,1},u_{2,2},u_{2,3},\ldots, u_{2,n},
         \cdots,
         u_{p_1,1},u_{p_1,2},u_{p_1,3},\ldots, u_{p_1,n},\\
        & v_{1,1},v_{1,2},v_{1,3},\ldots, v_{1,n},
          v_{2,1},v_{2,2},v_{2,3},\ldots, v_{2,n},
          \cdots,
          v_{q_1,1},v_{q_1,2},v_{q_1,3},\ldots, v_{q_1,n},\\
        & d_{1,1},d_{1,2},d_{1,3},\ldots,d_{1,p_{2}}\bigr) \in \mathbb{R}^{N},
\end{split}
\end{equation}
 where $$v_{k,1}=1, \text{ and } v_{k,l}=-1, \text{ for } l=2,3,\ldots, n,$$ and  the remaining coordinates are all zero.
 Then, $$D(S_{G_{1}}(n,n) \dot{\vee} G_{2})h^{(k,l)}=-2h^{(k,l)}, k=1,2,3,\ldots,q_{1},$$
  which yields a further collection of $q_{1}(n-1)$ linearly independent eigenvectors corresponding to the eigenvalue $-2.$
As the vectors $z^{(s,t)}$ and $h^{(k,l)}$ have disjoint supports, the two families together form a linearly independent set. Hence, $-2$ is an eigenvalue of $D(S_{G_{1}}(n,n) \dot{\vee} G_{2})$ with multiplicity at least $(n-1)(p_{1}+q_{1}).$\\

Let $X^{(j)}$ be an eigenvector of $A(G_{1})$ corresponding to the eigenvalue $\lambda_{1j}\neq r_{1}.$ Then, by the Perron-Frobenius theorem, $X^{(j)}$ is orthogonal to $J.$ By Lemma \ref{thelinegraph}, we have
\begin{align*}
    R(G_{1})R(G_{1})^{T}=&A(G_{1})+r_{1}I.\\
    R(G_{1})R(G_{1})^{T}X^{(j)}=&\big(A(G_{1})+r_{1}I\big)X^{(j)}\\
                    =&(\lambda_{1j}+r_{1})X^{(j)}.\\
               \text{ Also, }     B(G_{1})=&R(G_{1})^{T}R(G_{1})-2I.\\
                    B(G_{1})R(G_{1})^{T}X^{(j)}=&(R(G_{1})^{T}R(G_{1})-2I)R(G_{1})^{T}X^{(j)}\\
                    =&R(G_{1})^{T}(A(G_{1})+r_{1}I)X^{(j)}-2R(G_{1})^{T}X^{(j)}\\
                    =&(\lambda_{1j}+r_{1}-2)R(G_{1})^{T}X^{(j)}.
\end{align*}
Therefore, $R(G_{1})^{T}X^{(j)}$ is an eigenvector of $B(G_{1})$ with an eigenvalue $\lambda_{1j}+r_{1}-2,$ which is different from its regularity as $\lambda_{1j}\neq r_{1}.$ Then, again by the Perron-Frobenius theorem, $R(G_{1})X^{(j)}$ is orthogonal to $J.$
Consider $\psi_{j}=\begin{pmatrix}
    X^{(j)}_{p_{1}\times 1}\otimes J_{n\times 1}\\
    (R(G_{1})^{T}X^{(j)})_{q_{1}\times 1}\otimes J_{n\times 1}\\
    O_{p_{2}\times 1}
\end{pmatrix}.$ Then,
\begin{align*}
 & D(S_{G_{1}}(n,n) \dot{\vee} G_{2})\psi_{j}\\=&\begin{pmatrix} \big( 2(J-I)_{p_{1}\times p_{1}}\otimes J_{n \times n}\big)\big(X^{(j)}_{p_{1}\times 1}\otimes J_{n \times 1}\big)+\\
 2\widetilde{I}_{p_{1}n\times p_{1}n}\big(X^{(j)}_{p_{1}\times 1}\otimes J_{n\times 1}\big)+\\
 \Big( (3J-2R(G_{1}))_{p_{1}\times q_{1}}\otimes J_{n\times n}\Big) \Big( (R(G_{1})^{T}X^{(j)})_{q_{1}\times 1}\otimes J_{n\times 1} \Big)  \\
 \\
  \big( (3J-2R(G_{1})^{T})_{q_{1}\times p_{1}}X^{(j)}_{p_{1}\times 1}\big)\otimes (J_{n\times n} J_{n\times 1})+\\
  \big( (4(J-I)-2B(G_{1}))_{q_{1}\times q_{1}}\otimes J_{n\times n}\big) \big( (R(G_{1})^{T}X^{(j)})_{q_{1}\times 1}\otimes J_{n\times 1}\big)+\\
  2\widetilde{I}_{q_{1}n\times q_{1}n}\big( (R(G_{1})^{T}X^{(j)})_{q_{1}\times 1}\otimes J_{n\times 1}\Big)\\
  \\
  \big( J_{p_{2}\times p_{1}}\otimes J_{1\times n}\big) \big( X^{(j)}_{p_{1}\times 1}\otimes J_{n\times 1}\big)+\\
  2\big( J_{p_{2}\times q_{1}}\otimes J_{1\times n}\big) \big( (R(G_{1})^{T}X^{(j)})_{q_{1}\times 1} \otimes J_{n\times 1}\big)
  \end{pmatrix}\\
  =& \begin{pmatrix}\big( 2(J-I)_{p_{1}\times p_{1}}X^{(j)}_{p_{1}\times 1}\big) \otimes \big( J_{n\times n}J_{n\times 1}\big)+\\
  \big( (3J-2R)_{p_{1}\times q_{1}}(R^{T}X^{(j)})_{q_{1}\times 1}\big)\otimes (J_{n\times n}J_{n\times 1})+\\2\widetilde{I}_{p_{1}n\times p_{1}n}(X^{(j)}_{p_{1}\times 1}\otimes J_{n\times 1})\\
  \\
    \big( (3J-2R(G_{1})^{T})q_{1}\times p_{1}X^{(j)}_{p_{(i)}{1}\times 1}\big)\otimes\big(n J_{n\times 1} \big)+ \\
    \big( (4(J-I)-2B(G_{1}))_{q_{1}\times q_{1}} (R(G_{1})^{T}X^{(j)})_{q_{1}\times 1}\big)\otimes \big( nJ_{n\times 1}\big)+\\
    \big( I_{q_{1}\times q_{1}}\otimes 2(J-I)_{n\times n}\big) \big( (R^{T}X^{(j)})_{q_{1}\times 1}\otimes J_{n\times 1}\big)\\
    \\
    \big(J_{p_{2}\times p_{1}} X^{(j)}_{p_{1}\times 1}\big)\otimes \big( J_{1\times n}J_{n\times 1}\big)+\\
    \big( J_{p_{2}\times q_{1}}(R(G_{1})^{T}X^{(j)})_{q_{1}\times 1}\big)\otimes \big( J_{1\times n}J_{n\times 1}\big)\\
    \end{pmatrix}\\
    =&\begin{pmatrix}
        (-2X^{(j)}_{p_{1}\times 1})\otimes (nJ_{n\times 1})+(-2R(G_{1})R(G_{1})^{T})_{p_{1}\times 1}\otimes (nJ_{n\times 1}+\\
        \big( I_{p_{1}\times p_{1}}\otimes 2(J-I)_{n\otimes n}\big) \big( X^{(j)}_{p_{1}\times 1}\otimes J_{n\times 1}\big)\\
        \\
        \big( (3J-2R(G_{1})^{T})_{q_{1}\times p_{1}}X^{(j)}_{p_{1}\times 1}+\\
        (4(J-I)-2B(G_{1}))_{q_{1}\times q_{1}}(R(G_{1})^{T}X^{(j)})_{q_{1}\times 1}\big)\otimes (nJ_{n\times 1})+\\
        \big( (R(G_{1})^{T}X^{(j)})_{q_{1}\times 1}\otimes 2(n-1)J_{n\times 1}\\
        \\
       \big( O_{p_{2}\times 1}\otimes n\big)+2\big(O_{p_{2}\times 1}\otimes n\big)
    \end{pmatrix}\\
    =& \begin{pmatrix}
        \big( -2X^{(j)}_{p_{1}\times 1}-2(R(G_{1})R(G_{1})^{T}X^{(j)})_{p_{1}\times 1}\big) \otimes (nJ_{n\times 1})+\big( X^{(j)}_{p_{1}\times 1}\otimes 2(J-I)_{n\times n}J_{n\times 1}\big)\\
        \\
        \big( -2(R(G_{1})^{T}X^{(j)})_{q_{1}\times 1}-4(R(G_{1})^{T}X^{(j)})_{q_{1}\times 1}-2(\lambda_{1j}+r_{1}-2)(R(G_{1})^{T}X^{(j)})_{q_{1}\times 1}\big)\otimes nJ_{n\times 1}+\\2(n-1)\big((R(G_{1})^{T}X^{(j)})_{q_{1}\times 1}\otimes J_{n \times 1}\big)\\
        \\
        O_{p_{2}\times 1}
    \end{pmatrix}\\
    =&\begin{pmatrix}
        \big( -2X^{(j)}_{p_{1}\times 1}-2(A(G_{1})+r_{1}I)X^{(j)})_{p_{1}\times 1}\big)\otimes (nJ_{n\times 1})+2(n-1)(X^{(j)}_{p_{1}\times 1}\otimes J_{n\times 1})\\
       n(-2-2\lambda_{1j}-2r_{1}) \big((R(G_{1})^{T}X^{(j)})_{q_{1}\times 1}\otimes J_{n \times 1}\big)\\
        O_{p_{2}\times 1}
    \end{pmatrix}\\
    =&\begin{pmatrix}
        (-2-2\lambda_{1j}-2r_{1})(X^{(j)}_{p_{1}\times 1}\otimes nJ_{n\times 1})+2(n-1)(X^{(j)}_{p_{1}\times 1}\otimes J_{n\times 1})\\
        -2(n\lambda_{1j}+nr_{1}+1)\big((R(G_{1})^{T}X^{(j)})_{q_{1}\times 1}\otimes J_{n \times 1}\big)\\
        O_{p_{2}\times 1}
    \end{pmatrix}\\
    =&-2(n\lambda_{1j}+nr_{1}+1)\psi_{j}.
     \end{align*}
Thus, $\psi_{j}$ is an eigenvector of  $D(S_{G_{1}}(n,n) \dot{\vee} G_{2})$ with an eigenvalue $-2(n\lambda_{1j}+nr_{1}+1).$\\

Now, $-2$ is an eigenvalue of $B(G_{1})$ with multiplicity $q_{1}-p_{1}.$
Let $Z^{(i)}$ be an eigenvector of $B(G_{1})$ corresponding to the eigenvalue $-2.$ Then, $Z^{(i)}$ is orthogonal to $J,$ and by Lemma \ref{Beigenvalue}, $R(G_{1})Z^{(i)}=0.$
Now, consider the vector $\chi_{i}=\begin{pmatrix}
    O_{p_{1}\times 1}\otimes J_{n\times 1}\\
    Z^{(i)}_{q_{1}\times 1}\otimes J_{n\times 1}\\
    O_{p_{2}\times 1}
\end{pmatrix}.$
Then, \begin{align*}
   D(S_{G_{1}}(n,n) \dot{\vee} G_{2})\chi_{i}=&\begin{pmatrix}
  \big((3J-2R(G_{1}))_{p_{1}\times q_{1}}Z^{(i)}_{q_{1}\times 1}\big)\otimes J_{n\times n}J_{n\times 1}\\
  \big( (4(J-I)-2B)_{q_{1}\times q_{1}}\otimes J_{n\times n}\big)\big(Z^{(i)}_{q_{1}\times 1}\otimes J_{n\times 1} \big)+2\widetilde{I}_{q_{1}n \times q_{1}n}  \big(Z^{(i)}_{q_{1}\times 1}\otimes J_{n\times 1} \big)   \\
  \big( 2J_{p_{2}\times q_{1}}\otimes J_{1\times n}\big) \big(Z^{(i)}_{q_{1}\times 1} \otimes J_{n\times 1}\big) 
  \end{pmatrix}\\
  =&\begin{pmatrix}
  O_{p_{1}\times 1}\otimes nJ_{n\times 1}\\
  \big( (4(J-I)-2B )_{q_{1}\times q_{1}} Z^{(i)}_{q_{1}\times 1 }\big)\otimes J_{n\times n} J_{n\times 1}+2\widetilde{I}_{q_{1}n\times q_{1}n}(Z^{(i)}_{q_{1}\times 1}\otimes J_{n\times 1})\\
  \big(2J_{p_{2}\times q_{1}}Z^{(i)}_{q_{1}\times 1}\big) \otimes J_{1\times n}J_{n \times 1}
  \end{pmatrix}\\
  =&\begin{pmatrix}
  O_{p_{1}n\times 1}\\
  O_{q_{1}\times 1}\otimes nJ_{n\times 1}+\big(Z^{(i)}_{q_{1}\times 1}\otimes2(J-I)_{n \times n}J_{n \times 1}\big)\\
  O_{p_{2}\times 1}
  \end{pmatrix}\\
=&\begin{pmatrix}
  O_{p_{1}n\times 1}\\
  O_{q_{1}n\times 1}+2\big( Z^{(i)}_{q_{1}\times 1}\otimes (n-1)J_{n\times 1} \big)\\
  O_{p_{2}\times 1}
  \end{pmatrix}\\
=&\begin{pmatrix}
  O_{p_{1}n\times 1}\\
  O_{q_{1}n\times 1}+2(n-1)\big( Z^{(i)}_{q_{1}\times 1}\otimes J_{n\times 1}\big)\\
  O_{p_{2}\times 1}
  \end{pmatrix}\\
  =&2(n-1)\chi_{i}.
  \end{align*}
 Thus, $2(n-1)$ is an eigenvalue of $D(S_{G_{1}}(n,n) \dot{\vee} G_{2})$ with multiplicity at least $q_{1}-p_{1}.$
 
Now, let $Y^{(j)}$ be an eigenvector of $B(G_{1})$ corresponding to the eigenvalue $\lambda_{1j}+r_{2}-2,$  different from the regularity of $L(G_{1}).$ Then, $Y^{(j)}$ is orthogonal to $J$ and by Lemma \ref{thelinegraph}, $JR(G_{1})=2J.$\\
Let $\phi_{j}=\begin{pmatrix}
    (R(G_{1})Y^{(j)}_{p_{1}\times q_{1}}\otimes J_{n\times 1})_{p_{1}n\times 1}\\
    (-Y^{(j)}\otimes J_{n\times 1})_{q_{1}n\times 1}\\
    O_{p_{2}\times 1}
\end{pmatrix},$ then
\begin{align*}
    D(S_{G_{1}}(n,n) \dot{\vee} G_{2})\phi_{j}=&\begin{pmatrix}
        \big( 2(J-I)_{p_{1}\times p_{1}}\otimes J_{n\times n}\big)\big(RY^{(j)}_{p_{1}\times 1}\otimes J_{n \times 1}\big)+\\
        \big( 3J-2R)_{p_{1}\times q_{1}}(-Y^{(j)})_{q_{1}\times 1}\big)\otimes (J_{n\times n}J_{n\times 1})+\\
        2\widetilde{I}_{p_{1}n\times p_{1}n}\big(RY^{(j)}_{p_{1}\times 1}\otimes J_{n\times 1}\big)\\
        \\
          \big((3J-2R^{T})_{q_{1}\times p_{1}} \otimes J_{n\times n}\big)\big(RY^{(j)}_{p_{1}\times 1}\otimes J_{n\times 1}\big)+\\
          \Big( \big(4J-I)-2B(G_{1})\big)_{q_{1}\times q_{1}}\otimes J_{n\times n}\Big)\Big(-Y^{(j)}_{q_{1}\times 1}\otimes J_{n\times 1}\Big)+\\
          2\widetilde{I}_{q_{1}n\times q_{1}n}\big(-Y^{(j)}_{q_{1}\times 1}\otimes J_{n\times 1}\big)\\
          \\
         \big(J_{p_{2}\times p_{1}}\otimes J_{1\times n}\big)\big((R(G_{1})Y^{(j)})_{p_{1}\times 1}\otimes J_{n\times 1}\big)+\\
         2\big(J_{p_{2}\times q_{1}}\otimes J_{1\times n}\big)\big(-Y^{(j)}_{q_{1}\times 1}\otimes J_{n\times 1}\big)
         \end{pmatrix}\\
         =&\begin{pmatrix}
             \Big( \big(2(J-I)_{p_{1}\times p_{1}}RY^{(j)}_{p_{1}\times 1}\big)+\big( (3J-2R)_{p_{1}\times q_{1}}(-Y^{(j)})_{q_{1}\times 1}\big) \Big) \otimes \big(J_{n\times n}J_{n\times 1}\big)+\\
             \big( RY^{(j)}_{p_{1}\times 1}\otimes 2(J-I)_{n\times n}J_{n\times 1}\big)\\
             \\
             \Big( (3J-2R(G_{1})^{T})_{q_{1}\times p_{1}}(R(G_{1})Y^{(j)})_{p_{1}\times 1}\Big)\otimes \Big(J_{n\times n}J_{n\times 1}\Big)+\\
             \Big( \big( 4(J-I)-2B(G_{1})\big)_{q_{1}\times q_{1}}(-Y^{(j)})_{q_{1}\times 1}\Big)\otimes \Big( J_{n \times n } J_{n\times 1}\Big)+\\
             \Big(  I_{q_{1}\times q_{1}}\otimes 2(J-I)_{n\times n}\Big)\Big(-Y^{(j)}_{q_{1}\times 1}\otimes J_{n\times 1}\Big)\\
             \\
        \big( J_{p_{2}\times p_{1}} (R(G_{1})Y^{(j)})_{p_{1}\times 1}\big)\otimes \big(J_{1\times n}J_{n\times 1})+\\
        2\big(-J_{p_{2}\times q_{1}}Y^{(j)}_{q_{1}\times 1}\big)\otimes \big(J_{1\times n}J_{n\times 1}\big)
         \end{pmatrix}\\
         =&\begin{pmatrix}
             2(n-1)\big(RY^{(j)}_{p_{1}\times 1}\otimes J_{n\times 1}\big)\\
             2(n-1)\big(-Y^{(j)}_{q_{1}\times 1}\otimes J_{n\times 1}\big)\\
             O_{p_{2}\times 1}
         \end{pmatrix}\\
         =&2(n-1)\phi_{j}.
        \end{align*}
Thus, $2(n-1)$ is an eigenvalue of  $D(S_{G_{1}}(n,n) \dot{\vee} G_{2})$ with multiplicity at least $(q_{1}-p_{1})+(p_{1}-1)=q_{1}-1.$  \\

Now, let $\lambda_{2m}\neq r_{2}$ (for $m=2,3,\ldots,p_{2}$) be an eigenvalue of $A(G_{2})$ with an eigenvector $W^{(m)}$, and let $\eta_{m}=\begin{pmatrix}
    O_{p_{1}n\times 1}\\
    O_{q_{1}n\times 1}\\
    W^{(m)}
\end{pmatrix}.$ Then, $ D(S_{G_{1}}(n,n) \dot{\vee} G_{2})\eta_{m}=-(\lambda_{2m}+2)\eta_{m}.$ 
Thus,  $\eta_{m}$ is an eigenvector of $ D(S_{G_{1}}(n,n) \dot{\vee} G_{2})$ with an eigenvalue $-(\lambda_{2m}+2).$\\

Thus, we have obtained $(n-1)(p_{1}+q_{1})+(p_{1}-1)+(q_{1}-1)+(p_{2}-1)=np_{1}+nq_{1}+p_{2}-3$ eigenvalues. The remaining three eigenvalues are those
of the equitable quotient matrix of $D(S_{G_{1}}(n,n) \dot{\vee} G_{2})$,  $$M=\begin{pmatrix}
    2(np_{1}-1)&n(3q_{1}-2r_{1})&p_{2}\\
    n(3p_{1}-4)&4(q_{1}-r_{1})n+2(n-1)&2p_{2}\\
    p_{1}n&2q_{1}n&2(n-1)-r_{2}
\end{pmatrix}.$$ 
\end{proof}

The following lemma establishes the distance between any two vertices of $S_{G_{1}}(n,n) \veebar G_{2}.$
\begin{proposition}
    Let $V(G)=\{v_{1},v_{2},\ldots,v_{p}\}$ and $E(G)=\{e_{1},e_{2},\ldots, e_{q}\}.$ Then, the distance between any two distinct vertices, say $x$ and $y$, in $S_{G_{1}}(n,n) \veebar G_{2}$ is as follows,
\begin{equation*}
   d(x,y)=\begin{cases}
            2 & \text{ if } x=(v_{i},j),y=(v_{l},s)\text{ and }v_{i}\text{ is adjacent to } v_{l} \text{ in }G_{1}, \\
           4 & \text{ if } x=(v_{i},j),y=(v_{l},s)\text{ and }v_{i}\text{ is not adjacent to } v_{l} \text{ in }G_{1}, \\
         1 &\text{ if } x= (v_{i},j), y=(e_{l},s) \text{ and } v_{i} \text{ is incident with } e_{l} \text{ in } G_{1}, \\
             3 &\text{ if } x= (v_{i},j), y=(e_{l},s) \text{ and } v_{i} \text{ is not incident with } e_{l} \text{ in } G_{1}, \\
         2 & \text{ if }x= (v_{i},j) \text{ and }y \in V(G_{2}),\\
         2   &\text{ if } x= (e_{i},j), y=(e_{l},s),\\
         1 &\text{ if } x= (e_{i},j)\text { and } y \in V(G_{2}),\\
            1  &\text{ if }x,y\in V(G_{2}) \text{ and }x,y \text{ are adjacent in }V(G_{2}),\\
              2  &\text{ if }x,y\in V(G_{2}) \text{ and }x,y \text{ are not adjacent in }V(G_{2}).\\ 
            \end{cases}
\end{equation*}
\end{proposition}
\begin{lemma}\label{distance spectrum of generalized subdivision edge join}
 Let $G_{j}$ be an $r_{j}$-regular graph on $p_{j}$ vertices and $q_{j}$ edges with the adjacency matrix $A(G_{j})$ and adjacency spectrum $\{r_{j},\lambda_{j2},\lambda_{j3},\ldots,\lambda_{jp_{j}}\},\quad j=1,2.$ Then, the distance spectrum of $S_{G_{1}}(n,n) \veebar G_{2}$ consists of
  $-2$ with multiplicity $(n-2)p_{1}+nq_{1};$
     $-2-\lambda_{2m},m=2,3,\ldots,p_{2};$
$-2(nt+1),t=\frac{(\lambda_{1i}+1)\pm \sqrt{(\lambda_{1i}+1)^2+4(\lambda_{1i}+r_{1})}}{2},i=2,3,\ldots, p_{1};$
 and the eigenvalues of the matrix, $$\begin{pmatrix}
         (4(p_{1}-1)-2r_{1})n+2(n-1)&(3q_{1}-2r_{1})n&2p_{2}\\
         n(3p_{1}-4)&2(q_{1}n-1)&p_{2}\\
         2p_{1}n&q_{1}n&2(p_{2}-1)-r_{2}
     \end{pmatrix}.$$
\end{lemma}
\begin{proof} 
By suitable labeling of the vertices of $(S_{G_{1}}(n,n) \veebar G_{2}),$ the distance matrix $D(S_{G_{1}}(n,n) \veebar G_{2})$ can be represented as
    $$\begin{pmatrix}
\big(4(J-I)-2A(G_{1})\big)_{p_{1}\times p_{1}}\otimes J_{n\times n}+\widetilde{I}_{p_{1}n\times p_{1}n}&(3J-2R(G_{1}))_{p_{1}\times q_{1}}\otimes J_{n\times n}&2J_{p_{1}n\times p_{2}}\\
(3J-2R(G_{1})^{T})_{q_{1}\times p_{1}}\otimes J_{n\times n}&\big( 2(J-I)_{q_{1}\times q_{1}}\otimes J_{n\times n}\big)+\widetilde{I}_{q_{1}n\times q_{1}n} &J_{q_{1}n\times p_{2}}\\
2J_{p_{2}\times p_{1}n}&2J_{p_{2}\times q_{1}n}&\big(2(J-I)-A(G_{2})\big)_{p_{2}\times p_{2}}
 \end{pmatrix}.$$
Let $\lambda_{1i}\neq r_{1}$ be an eigenvalue of $A(G_{1})$ with an eigenvector $X^{(i)}.$ Then, by Perron-Frobenius Theorem, $X^{(i)}$ is orthogonal to $J.$ Consider $\phi_{i}=\begin{pmatrix}
    tX^{(i)}_{p_{1}\times 1}\otimes J_{n\times 1}\\
    \\
    R(G_{1})^{T}_{q_{1}\times 1}X^{(i)}_{p_{1}\times 1}\\
    \\
    O_{p_{2}\times 1}
\end{pmatrix}$
is an eigenvector of $D(S_{G_{1}}(n,n) \veebar G_{2}).$ If $\mu_{i}$ is an eigenvalue of $D(S_{G_{1}}(n,n) \veebar G_{2})$
with eigenvector $\phi_{i},$ then
$D(S_{G_{1}}(n,n) \veebar G_{2}) \phi_{i}=\mu_{i}\phi_{i}.$
That is, 
\begin{align*}
  \begin{pmatrix}
      \Big( \big( 4(J-I)-2A(G_{1})\big)_{p_{1}\times p_{1}}\otimes J_{n\times n}\Big)\Big( tX_{p_{1}\times 1}^{(i)}\otimes J_{n\times 1}\Big)+\\
      \widetilde{I}_{p_{1}n\times p_{1}n}\big( tX_{p_{1}\times 1}^{(i)}\otimes J_{n\times 1}\big)+\\ 
      \big( (3J-2R(G_{1})_{p_{1}\times q_{1}}(R(G_{1})^{T}X^{(i)})_{q_{1}\times 1}\big)\otimes \big(J_{n\times n}J_{n\times 1}\big)\\
      \\
      \Big( \big(3J-2R(G_{1})\big)_{q_{1}\times p_{1}}(tX^{(i)})_{p_{1}\times 1}\Big)\otimes \Big(J_{n\times n}J_{n\times 1}\Big)+\\
      \Big( 2(J-I)_{q_{1}\times q_{1}}\otimes J_{n\times n}\Big)\Big((R(G_{1})^{T}X^{(i)})_{q_{1}\times 1}\otimes J_{n\times 1}\Big)\\
      \widetilde{I}_{q_{1}n\times q_{1}n}\Big( (R(G_{1})^{T}X^{(i)})_{q_{1}\times 1}\otimes J_{n\times 1}\Big)
      \\
      \\
      \Big(2J_{p_{2}\times p_{1}}\otimes J_{1\times n}\Big)\Big(tX_{p_{1}\times 1}^{(i)}\otimes J_{n\times 1}\Big)+
      \\
      \Big( J_{p_{2}\times q_{1}}\otimes J_{1}\times 1\Big)\Big( (R(G_{1})^{T}X^{(i)})_{q_{1}\times 1}\otimes J_{n\times 1}\Big)
    \end{pmatrix}
    =\mu_{i}\begin{pmatrix}
    tX^{(i)}_{p_{1}\times 1}\otimes J_{n\times 1}\\
    \\
    R(G_{1})^{T}_{q_{1}\times 1}X^{(i)}_{p_{1}\times 1}\\
    \\
    O_{p_{2}\times 1}
\end{pmatrix}\\
\begin{pmatrix}
\Big( \big(4(J-I)-2A(G_{1})\big)_{p_{1}\times p_{1}}tX_{p_{1}\times 1}^{(i)}\Big)\otimes \Big( J_{n\times n}J_{n\times 1}\Big)+\\\Big( \big(3J-2R(G_{1})\big)_{p_{1}\times q_{1}}(R(G_{1})^{T}X^{(i)})_{q_{1}\times 1}\Big)\otimes \Big( J_{n\times n}J_{n\times 1}\Big)+\\\Big( \widetilde{I}_{p_{1}\times p_{1}}\otimes 2(J-I)_{n\times n}\Big)\Big(tX_{p_{1}n\times 1}^{(i)}\otimes J_{n\times 1}\Big)\\
\\
\Big( -2t(R(G_{1})^{T}X^{(i)})_{q_{1}\times 1}\Big)\otimes \Big(nJ_{n\times 1}\Big)+\Big(2(J-I)R(G_{1})^{T}_{q_{1}\times 1}\Big)\otimes \Big(nJ_{n\times 1}\\
\Big((R(G_{1})^{T}X^{(i)})_{q_{1}\times 1}\otimes 2(J-I)_{n\times n}J_{n\times 1}
\\
\\
2t\Big((JX^{(i)})_{p_{2}\times 1}\otimes nJ_{1\times 1}\Big)+\\
\Big(J_{p_{2}\times q_{1}}(R(G_{1})^{T}X^{(i)})_{q_{1}\times 1}\Big)\otimes \big( J_{1\times n}J_{n\times 1}\Big)
\end{pmatrix}
 =\mu_{i}\begin{pmatrix}
    tX^{(i)}_{p_{1}\times 1}\otimes J_{n\times 1}\\
    \\
    R(G_{1})^{T}_{q_{1}\times 1}X^{(i)}_{p_{1}\times 1}\\
    \\
    O_{p_{2}\times 1}
\end{pmatrix}&\\
\begin{pmatrix}
\Big( t\big(-2n-2-2\lambda_{1i}n\big)-2n\big(\lambda_{1i}+r_{1}\big)\Big)\Big(X_{p_{1}\times 1}^{(i)}\otimes J_{n\times 1}\Big)\\
\\
 (-2nt-2)\Big(R(G_{1})^{T}_{q_{1}\times 1}X^{(i)}_{p_{1}\times 1}\Big)\\
 \\
 O_{p_{2}\times 1}
 \end{pmatrix}
= \mu_{i}\begin{pmatrix}
    tX^{(i)}_{p_{1}\times 1}\otimes J_{n\times 1}\\
    \\
    R(G_{1})^{T}_{q_{1}\times 1}X^{(i)}_{p_{1}\times 1}\\
    \\
    O_{p_{2}\times 1}
\end{pmatrix}.\\
\end{align*}
Thus, 
\begin{align*}
    t(-2n-2-2n\lambda_{1i})-2n(\lambda_{1i}+r_{1})=&\mu_{i}t\\
    -2nt-2=&\mu_{i}.
\end{align*}
Hence, \begin{align*}
    t^{2}-t(\lambda_{1i}+1)-(\lambda_{1i}+r_{1})=0.\\
    \text{ Therefore, $t$ has two values, } 
    t_{1}=&\frac{(\lambda_{1i}+1)+ \sqrt{(\lambda_{1i}+1)^2+4(\lambda_{1i}+r_{1})}}{2},\\
    t_{2}=&\frac{(\lambda_{1i}+1)- \sqrt{(\lambda_{1i}+1)^2+4(\lambda_{1i}+r_{1})}}{2}.
\end{align*}
That is, we get two values of $\mu_{i},$ $-2(nt_{1}+1)$ and $-2(nt_{2}+1),$ corresponding to each eigenvalue $\lambda_{1i}\neq r_{1}$ of $G_{1}.$ Thus, we obtain a  total of $2p_{1}-2$ eigenvalues of $D(S_{G_{1}}(n,n) \veebar G_{2}).$

Now, consider the vector (\ref{eq:eigenvector_gsvj one}). Then we get, 
$$D(S_{G_{1}}(n,n) \veebar G_{2})z^{(s,t)}=-2z^{(s,t)}, s=1,2,3,\ldots,p_{1}, t=2,3,\ldots n. $$ 
This yields $(n-1)p_{1}$ linearly independent eigenvectors of $D(S_{G_{1}}(n,n) \veebar G_{2})$  associated with the eigenvalue $-2.$
Similarly, for the vector (\ref{eq:eigenvector_gsvj two}), we have
$$D(S_{G_{1}}(n,n) \veebar G_{2})h^{(k,l)}=-2h^{(k,l)}, k=1,2,\ldots,q_{1},l=2,3,\ldots n$$
 Thereby, we obtain a further collection of $q_{1}(n-1)$ linearly independent eigenvectors corresponding to the eigenvalue $-2.$\\
Let $Z_{q_{1}\times 1}^{(d)}$ be an eigenvector of $B(G_{1})$ corresponding to the eigenvalue $-2,$ with multiplicity $q_{1}-p_{1}.$ By Lemma \ref{Beigenvalue}, $R(G_{1})_{p_{1}\times q_{1}}Z_{q_{1}\times 1}^{(d)}=O.$
Let $\psi_{d}=\begin{pmatrix}
    O_{p_{1}\times 1}\\
    Z_{q_{1}\times 1}^{(d)}\otimes J_{n\times 1}\\
    O_{p_{2}\times 1}
\end{pmatrix}.$
Then, \begin{align*}
D(S_{G_{1}}(n,n) \veebar G_{2})\psi_{d}=&\begin{pmatrix}
    \big( (3J-2R(G_{1}))_{p_{1}\times q_{1}} Z_{q_{1}\times 1}^{(d)}\big)\otimes \big( J_{n\times n}J_{n\times 1}\big)\\
    \big(2(J-I)_{q_{1}\times q_{1}}\otimes J_{n\times n}\big) \big( Z_{q_{1}\times 1}^{(d)}\otimes J_{n\times 1}\big)+2\widetilde{I}_{q_{1}n\times q_{1}n}\big( Z_{q_{1}\times 1}^{(d)}\otimes J_{n\times 1}\big)                \\
    \big( J_{p_{2}\times q_{1}}\otimes J_{1\times n}\big)\big(Z_{q_{1}\times 1}^{(d)}\otimes J_{n\times 1}\big)
\end{pmatrix}\\
=&\begin{pmatrix}
    O_{p_{1}\times 1}\otimes nJ_{n\times 1}\\
    \big( 2(J-I)_{q_{1}\times q_{1}}Z_{q_{1}\times 1}^{(d)}\big)\otimes \big( J_{n\times n}J_{n\times 1}\big)+\big( I_{q_{1}\times q_{1}}\otimes 2(J-I)_{n\times n}\big)\big(Z_{q_{1}\times 1}^{(d)}\otimes J_{n\times 1}\big)\\
    \big( J_{p_{2}\times q_{1}}Z_{q_{1}\times 1}^{(d)}\big)\otimes \big(J_{1\times n}J_{n\times 1}\big)
\end{pmatrix}\\
=&\begin{pmatrix}
    O_{p_{1}\times 1}\otimes J_{n\times 1}\\
    \big( -2nZ_{q_{1}\times 1}^{(d)}\otimes J_{n\times 1}\big)+\big(Z_{q_{1}\times 1}^{(d)}\otimes (2(J-I)_{n\times n}J_{n\times 1})\big)\\
    O_{p_{2}\times 1}\otimes J_{1\times 1}
\end{pmatrix}\\
=&\begin{pmatrix}
    O_{p_{1}\times 1}\otimes J_{n\times 1}\\
    -2\big(Z_{q_{1}\times 1}^{(d)}\otimes J_{n\times 1}\big)\\
    O_{p_{2}\times 1}
\end{pmatrix}
=-2\psi_{d},
\end{align*}
which again yields  $q_{1}-p_{1}$ linearly independent eigenvectors corresponding to the eigenvalue $-2.$\\
Finally, the pairwise disjointness of the coordinate support of the vectors  $\{z^{(s,t)}\}, \{h^{(k,l)}\}$ and $\{\psi_{d}\}, \text{( for all } s,t,k,l,d),$ ensures that they form a linearly independent set. Hence, $-2$ is an eigenvalue of $D(S_{G_{1}}(n,n) \veebar G_{2})$ with multiplicity at least $(n-1)p_{1}+(n-1)q_{1}+(q_{1}-p_{1})=(n-2)p_{1}+nq_{1}.$\\

Now, let $\lambda_{2m}\neq r_{2}$ (for $m=2,3,\ldots, p_{2}$) be an eigenvalue of $A(G_{2})$ with an eigenvector $W^{(m)}_{p_{2}\times 1}.$ Let $\chi_{m}=\begin{pmatrix}
    O_{p_{1}n\times 1}\\
    O_{q_{1}n\times 1}\\
    W_{p_{2}\times 1}^{(m)}
\end{pmatrix}.$
Then, \begin{align*}
D(S_{G_{1}}(n,n) \veebar G_{2})\chi_{m}=&\begin{pmatrix}
    2J_{p_{1}n\times p_{2}}W^{(m)}_{p_{2}\times 1}\\
     J_{q_{1}n\times p_{2}}W^{(m)}_{p_{2}\times 1}\\
     \big( 2(J-I)-A(G_{2})\big)_{p_{2}\times p_{2}}W^{(m)}_{p_{2}\times 1}
\end{pmatrix}\\
=&\begin{pmatrix}
O_{p_{1}n\times 1}\\
O_{q_{1}n\times 1}\\
(-2-\lambda_{2m})W^{(m)}_{p_{2}\times 1}
\end{pmatrix}
=(-2-\lambda_{2m})\chi_{m}.
\end{align*}
This shows that $(-2-\lambda_{2m})$ is an eigenvalue of $D(S_{G_{1}}(n,n) \veebar G_{2})$ corresponding to the
eigenvalue $\lambda_{2m}$ of $A(G_{2}).$\\
Thus, we get $p_{1}n+q_{1}n+p_{2}-3$ eigenvalues.
The remaining $3$ eigenvalues are given by the eigenvalues of the equitable quotient matrix of $D(S_{G_{1}}(n,n) \veebar G_{2}),$
 $$M=\begin{pmatrix}
         (4(p_{1}-1)-2r_{1})n+2(n-1)&(3q_{1}-2r_{1})n&2p_{2}\\
         n(3p_{1}-4)&2(q_{1}n-1)&p_{2}\\
         2p_{1}n&q_{1}n&2(p_{2}-1)-r_{2}
     \end{pmatrix}.$$

\end{proof}

\section{Applications}
In this section, we construct some new families of $D$-equienergetic graphs and $D$-cospectral graphs of diameter $4$.\\

Lemmas \ref{distance spectrum of generalized subdivision vertex join} and \ref{distance spectrum of generalized subdivision edge join} demonstrate that the $D$-spectrum of
$S_{G_{1}}(n,n)\dot{\vee}G_{2}$ and $S_{G_{1}}(n,n)\veebar G_{2}$, 
are determined by the orders, regularity and the adjacency spectra of $G_{1}$ and $G_{2},$
together with the eigenvalue  $-2$ of $B(G_{1}),$ the adjacency matrix of the line graph of $G_{1}.$ In view of this, Lemmas \ref{distance spectrum of generalized subdivision vertex join} and \ref{distance spectrum of generalized subdivision edge join} yield the following theorems.
\begin{theorem}
Let $G_{1}$ and $G_{2}$ be regular cospectral graphs of the same order and regularity, and $H$ be any regular graph. Then, 
\begin{enumerate}
    \item[\rm(i)] $S_{G_{1}}(n,n)\dot{\vee}H$ and $S_{G_{2}}(n,n)\dot{\vee}H$ are $D$-cospectral.
    \item[\rm(ii)]  $S_{G_{1}}(n,n)\veebar H$ and $S_{G_{2}}(n,n)\veebar H$ are $D$-cospectral.
\end{enumerate}
\end{theorem}
\begin{theorem}\label{thm:distance equienergetic-both}
Let $G$ be an $r$-regular graph on $p$ vertices and $q$ edges. For a fixed $t\in\mathbb{N}$, let $\mathcal{T}$ be the family of all integer partitions $T=(t_1,t_2,\ldots,t_s)$ of $t$ into parts $t_i\geq 3$ $(1\leq i\leq s)$, and for each such $T\in\mathcal{T}$, let
$$C_T=C_{t_1}\cup C_{t_2}\cup\cdots\cup C_{t_s}$$
be the corresponding disjoint union of cycles. Then,
\begin{enumerate}
    \item[\rm(i)] the graphs $S_G(n,n)\,\dot{\vee}\,C_T,\ T\in\mathcal{T}$, form a family of $D$-equienergetic graphs,
    \item[\rm(ii)] the graphs $S_G(n,n)\,\veebar\,C_T,\ T\in\mathcal{T}$, form  a family of $D$-equienergetic graphs.
\end{enumerate}
\end{theorem}
\begin{proof}
       Let $\lambda_{1}=r_{1}\geq \lambda_{2} \geq \ldots \geq \lambda_{p}$ and   $\mu_{1}=2\geq \mu_{2} \geq \ldots \geq \mu_{t}$ be the adjacency eigenvalues of $G$ and $C_{T}$ respectively.
       \begin{enumerate}
           \item[\rm(i)]Then, by Lemma \ref{distance spectrum of generalized subdivision vertex join}, the $D$-eigenvalues of $S_{G}(n,n)\dot{\vee} C_{T}$ are, 
 $-2$ with multiplicity $(n-1)(p+q);$
      $2(n-1)$ with multiplicity $q-1;$
      $-2(n\lambda_{j}+nr+1), j=2,3,\ldots,p;$
      $-2-\mu_{i},i=2,3,\ldots t;$
      the eigenvalues of the matrix, $$\begin{pmatrix}
         2(np-1)&n(3q-2r)&t\\
         n(3p-4)&4(q-r)n+2(n-1)&2t\\
         pn&2qn&2(t-1)-2
        \end{pmatrix},$$
for every partition $T$ of $t.$
Since $\mu_{i}\geq -2,$ for $i=2,3,\ldots, t,$ we have, $-(\mu_{i}+2)\leq 0.$ Therefore, \begin{align*}
    \sum_{i=2}^{t}|-2-\mu_{i}|=& \sum_{i=2}^{t}(2+\mu_{i})\\
                              =&2(t-1)+ \sum_{i=2}^{t}\mu_{i}\\
                              =&2t-4.
\end{align*}
Thus, for every partitions $T_{i},T_{j}\in \mathcal{T,}$ 
$$E_{D}(S_{G}(n,n)\dot{\vee} C_{T_{i}})=E_{D}(S_{G}(n,n)\dot{\vee} C_{T_{j}}).$$
\item[\rm(ii)] By Lemma \ref{distance spectrum of generalized subdivision vertex join}, the $D$-eigenvalues of $S_{G}(n,n)\veebar  C_{T}$ are, 
$-2$ with multiplicity $(n-2)p+nq$; 
$-2-\mu_{i},i=2,3,\ldots,t$;
 $-2(nt_{j}+1), t=\frac{(\lambda_{i}+1)\pm \sqrt{(\lambda_{i}+1)^2+4(\lambda_{i}+r_{1})}}{2},j=1,2$;
and the eigenvalues of the matrix $$\begin{pmatrix}
         (4(p_{1}-1)-2r_{1})n+2(n-1)&(3q_{1}-2r_{1})n&2t\\
         n(3p_{1}-4)&2(q_{1}n-1)&t\\
         2p_{1}n&q_{1}n&2(t-1)-2
     \end{pmatrix},$$
for every partition $T$ of $t.$ Therefore, as in part $(i)$ we get, \begin{align*}
    \sum_{i=2}^{t}|-2-\mu_{i}|=& 2t-4.
\end{align*}
Hence, 
$$E_{D}(S_{G}(n,n)\veebar C_{T_{i}})=E_{D}(S_{G}(n,n)\veebar C_{T_{j}}),$$
 for every partitions $T_{i},T_{j}\in \mathcal{T}.$
       \end{enumerate}
\end{proof}

\newpage
Table \ref{table1} lists some families of $D$-equienergetic graphs constructed using the generalized subdivision vertex join and generalized subdivision edge join of graphs.

\begin{table}[h]
\centering
\begin{tabular}{ccccc}
\hline
 No & $t$  & n & $D$-equienergetic graphs &$D$-energy \\ 
\hline
$1$&$6$&$12$&$S_{K_{2}}(2,2)\dot{\vee}C_{6},\hspace{0.1cm}S_{K_{2}}(2,2)\dot{\vee}(C_{3}\cup C_{3})$&$37.1324$\\
$2$&$6$&$15$&$S_{K_{2}}(3,3)\dot{\vee}C_{6},\hspace{0.1cm}S_{K_{2}}(3,3)\dot{\vee}(C_{3}\cup C_{3})$&$49.4684$\\
$3$&$7$&$13$&$S_{K_{2}}(2,2)\dot{\vee}C_{7},\hspace{0.1cm}S_{K_{2}}(2,2)\dot{\vee}(C_{4}\cup C_{3})$&$40.9979$\\
$4$&$7$&$16$&$S_{K_{2}}(3,3)\dot{\vee}C_{7},\hspace{0.1cm}S_{K_{2}}(3,3)\dot{\vee}(C_{4}\cup C_{3})$&$53.3158$\\
$5$&$8$&$14$&$S_{K_{2}}(2,2)\dot{\vee}C_{8},\hspace{0.1cm}S_{K_{2}}(2,2)\dot{\vee}(C_{4}\cup C_{4}),\hspace{0.1cm} S_{K_{2}}(2,2)\dot{\vee}(C_{5}\cup C_{3})$&$44.8908$\\
$6$&$8$&$17$&$S_{K_{2}}(3,3)\dot{\vee}C_{8},\hspace{0.1cm}S_{K_{2}}(3,3)\dot{\vee}(C_{4}\cup C_{4}),\hspace{0.1cm} S_{K_{2}}(3,3)\dot{\vee}(C_{5}\cup C_{3})$&$57.1908$\\
$7$&$6$&$12$&$S_{K_{2}}(2,2)\veebar C_{6},\hspace{0.1cm} S_{K_{2}}(2,2)\veebar(C_{3}\cup C_{3})$&$37.7234$\\
$8$&$6$&$15$&$S_{K_{2}}(3,3)\veebar C_{6},\hspace{0.1cm} S_{K_{2}}(3,3)\veebar(C_{3}\cup C_{3})$&$50.1056$\\
$9$&$7$&$13$&$S_{K_{2}}(2,2)\veebar C_{7},\hspace{0.1cm} S_{K_{2}}(2,2)\veebar(C_{4}\cup C_{3})$&$41.5641$\\
$10$&$7$&$16$&$S_{K_{2}}(3,3)\veebar C_{7},\hspace{0.1cm} S_{K_{2}}(3,3)\veebar(C_{4}\cup C_{3})$&$53.9442$\\
$11$&$8$&$14$&$S_{K_{2}}(2,2)\veebar C_{8},\hspace{0.1cm} S_{K_{2}}(2,2)\veebar(C_{4}\cup C_{4}),\hspace{0.1cm} S_{K_{2}}(2,2)\veebar(C_{5}\cup C_{3})$&$45.4306$\\
$12$&$8$&$17$&$S_{K_{2}}(3,3)\veebar C_{8},\hspace{0.1cm} S_{K_{2}}(3,3)\veebar(C_{4}\cup C_{4}),\hspace{0.1cm}  S_{K_{2}}(3,3)\veebar(C_{5}\cup C_{3})$&$57.8046$\\
\hline
\end{tabular}
\caption{ Some $D$-equienergetic graphs  of diameter 4.}
\label{table1}
\end{table}

Lemmas  \ref{distance spectrum of generalized subdivision vertex join}  and \ref{distance spectrum of generalized subdivision edge join} determine the distance spectra of
$S_{G_1}(n,n)\dot{\vee}(G_2\cup G_3)$ and
$S_{G_1}(n,n)\veebar(G_2\cup G_3)$, respectively, when
$G_2$ and $G_3$ are regular graphs of the same degree.

The following theorems extend these results to the more general
case in which $G_2$ and $G_3$ are regular graphs of distinct
degrees.

\begin{theorem}
    For $i=1,2,3,$ let $G_{i}$ be an $r_{i}$-regular graph with  $p_{i}$ vertices and  $q_{i}$ edges whose adjacency spectrum is given by $\{r_{i},\lambda_{i1},\lambda_{i2},\ldots,\lambda_{ip_{i}}\}.$  If $r_{2}\neq r_{3},$ then the $D$-spectrum of $S_{G_{1}}(n,n) \dot{\vee}( G_{2}\cup G_{3})$ consists of 
         $-2$ with multiplicity $(n-1)(p_{1}+q_{1})$;
         $2(n-1)$ with multiplicity $q_{1}-1$;
           $-2-\lambda_{ij}, i=1,2 ;j=2,3,\ldots,p_{2}$;
     and the eigenvalues of the matrix, $$\begin{pmatrix}
         2(np_{1}-1)&n(3q_{1}-2r_{1})&p_{2}&p_{3}\\
         n(3p_{1}-4)&4(q_{1}-r_{1})n+2(n-1)&2p_{2}&2p_{3}\\
         p_{1}n&2q_{1}n&2(p_{2}-1)-r_{2}&2p_{3}\\
         p_{1}n&2q_{1}n&2p_{2}&2(p_{3}-1)-r_{3}
        \end{pmatrix}.$$
   
\end{theorem}
\begin{proof}
By an appropriate labeling of vertices of $S_{G_{1}}(n,n) \dot{\vee}( G_{2}\cup G_{3}),$ the distance matrix of $D(S_{G_{1}}(n,n) \dot{\vee}( G_{2}\cup G_{3}))$ can be represented as
 $$\resizebox{\textwidth}{!}{$\begin{pmatrix}
       \big(2(J-I)_{p_{1}\times p_{1}}\otimes J_{n\times n} \big)+2 \widetilde{I}_{p_{1}n\times p_{1}n}   & (3J-2R(G_{1}))_{p_{1}\times q_{1}} \otimes J_{n\times n}   &  J_{p_{1}n\times p_{2}}&J_{p_{1}n\times p_{3}}\\
       \big( (3J-2R(G_{1})^{T})_{q_{1}\times p_{1}}\otimes J_{n\times n}  \big)  & \big( (4(J-I)-2B(G_{1}))_{q_{1}\times q_{1}}\otimes J_{n\times n}\big)+2\widetilde{I}_{q_{1}n\times q_{1}n} & 2J_{q_{1}n\times p_{2}} &2J_{q_{1}n\times p_{3}}\\
       J_{p_{2}\times p_{1}n} &2J_{p_{2}\times q_{1}n} &\big( 2(J-I)-A(G_{2})\big)_{p_{2}\times p_{2}}&J_{p_{2}\times p_{3}}\\
       J_{p_{3}\times p_{1}n}&2J_{p_{3}\times q_{1}n}&2J_{p_{3}\times p_{2}}&(2(J-I)-A(G_{3}))_{p_{3}\times p_{3}}
   \end{pmatrix}$}.$$
    The proof is analogous to that of Lemma \ref{distance spectrum of generalized subdivision vertex join}.
\end{proof}

\begin{theorem}
    For $i=1,2,3,$ let $G_{i}$ be an $r_{i}$-regular graph with order $p_{i}$ and size $q_{i}$ whose adjacency spectrum is given by $\{r_{i},\lambda_{i1},\lambda_{i2},\ldots,\lambda_{ip_{i}}\}.$ If $r_{2}\neq r_{3},$ then the $D$-spectrum of $S_{G_{1}}(n,n) \veebar ( G_{2}\cup G_{3})$ consists of
     $-2$ with multiplicity $(n-2)p_{1}+nq_{1}$;
      $-2-\lambda_{ij},j=2,3,\ldots,p_{2},i=2,3$;
      $-2(nt+1),t=\frac{(\lambda_{1i}+1)\pm \sqrt{(\lambda_{1i}+1)^2+4(\lambda_{1i}+r_{1})}}{2},i=2,3,\ldots, p_{1} $;
     and  the eigenvalues of the matrix, $$\begin{pmatrix}
         (4(p_{1}-1)-2r_{1})n+2(n-1)&(3q_{1}-2r_{1})n&2p_{2}&2p_{3}\\
         n(3p_{1}-4)&2(q_{1}n-1)&p_{2}&p_{3}\\
         2p_{1}n&q_{1}n&2(p_{2}-1)-r_{2}&2p_{3}\\
         2p_{1}n&q_{1}n&2p_{2}&2(p_{3}-1)-r_{3}
     \end{pmatrix}.$$

\end{theorem}
\begin{proof}
By an appropriate labeling of vertices of $S_{G_{1}}(n,n) \veebar ( G_{2}\cup G_{3}),$ the distance matrix of $D(S_{G_{1}}(n,n) \dot{\vee}( G_{2}\cup G_{3}))$ can be represented as
    $$\scriptsize{\begin{pmatrix}
\big(4(J-I)-2A(G_{1})\big)_{p_{1}\times p_{1}}\otimes J_{n\times n}+\widetilde{I}_{p_{1}n\times p_{1}n}&(3J-2R(G_{1}))_{p_{1}\times q_{1}}\otimes J_{n\times n}&2J_{p_{1}n\times p_{2}}&2J_{p_{1}n\times p_{3}}\\
(3J-2R(G_{1})^{T})_{q_{1}\times p_{1}}\otimes J_{n\times n}&\big( 2(J-I)_{q_{1}\times q_{1}}\otimes J_{n\times n}\big)+\widetilde{I}_{q_{1}n\times q_{1}n} &J_{q_{1}n\times p_{2}}&J_{q_{1}n\times p_{3}}\\
2J_{p_{2}\times p_{1}n}&2J_{p_{2}\times q_{1}n}&\big(2(J-I)-A(G_{2})\big)_{p_{2}\times p_{2}}&2J_{p_{2}\times p_{3}}\\
2J_{p_{3}\times p_{1}n}&J_{p_{3}\times q_{1}n} &2J_{p_{3}\times p_{2}}&\big(2(J-I)-A(G_{3})\big)_{p_{3}\times p_{3}}
 \end{pmatrix}}.$$
  The proof follows from the same arguments as in Theorem \ref{distance spectrum of generalized subdivision edge join}.
\end{proof}

\begin{theorem}
 Let $G$ be an $r$-regular graph on  $p$ vertices and  $q$ edges, and $H$ be a regular graph with the smallest eigenvalue at least $-2$. For a fixed $t\in \mathbb{N},$ let $\mathcal{T}$ be the family of all integer partitions $T=(t_{1},t_{2},\ldots, t_{s})$ of $t$ into parts $t_{i}\geq 3,(1\leq i\leq s).$  For each such $T\in \mathcal{T},$ let $C_{T}=C_{t_{1}}\cup C_{t_{2}}\cup \cdots \cup C_{t_{s}}$ be the corresponding disjoint union of cycles. Then,   \begin{enumerate}
 \item[\rm(i)] the graphs $S_{G}(n,n)\dot{\vee} (C_{T}\cup H),T\in \mathcal{T}$ form a family of $D$-equienergetic graphs,
 \item[\rm(ii)] the graphs $S_{G}(n,n)\veebar  (C_{T}\cup H),T\in \mathcal{T}$ form a family of $D$-equienergetic graphs.
 \end{enumerate}
  
\end{theorem}
\section*{Acknowledgement}
The research of Anjitha Ashokan is supported by the University Grants Commission of India under the beneficiary code BININ05086971.
\section*{Declarations}
 On behalf of all authors, the corresponding author states that there is no conflict of interest.

\bibliographystyle{plain}
 \bibliography{reference}

@article{cvetkovic1980spectra,
    author= "Cvetkovi\'{c}, Drago\v{s} M. and Doob, Michael and Sachs,
              Horst",
     title= "{Spectra of graphs Theory and application.}",
 journal= "vol. 87, Academic Press, Inc,
              New York-London",
      year= "1980",
       MRCLASS = "05Cxx (05-01 92A40)",
  MRNUMBER = "572262",
MRREVIEWER = "N.\ Trinajsti\'{c}",
}

@article{zbMATH03363442,
 author = {Graham, R. L. and Pollak, H. O.},
 title = {On the addressing problem for loop switching},
 fjournal = {Bell System Technical Journal},
 journal = {Bell Syst. Tech. J.},
 issn = {0005-8580},
 volume = {50},
 pages = {2495--2519},
 year = {1971},
 language = {English},
 doi = {10.1002/j.1538-7305.1971.tb02618.x},
 zbMATH = {3363442},
 Zbl = {0228.94020}
}

@article{balasubramanian1990computer,
  title={Computer generation of distance polynomials of graphs},
  author={Balasubramanian, K},
  journal={J. Comput. Chem.},
  volume={11},
  number={7},
  pages={829--836},
  year={1990},
  publisher={Wiley Online Library}
}

@article{indulal2008distance,
  title={On distance energy of graphs},
  author={Indulal, Gopalapillai and Gutman, Ivan and Vijayakumar, Ambat},
  journal={MATCH Commun. Math. Comput. Chem},
  volume={60},
  number={2},
  pages={461--472},
  year={2008}
}

@article{zbMATH07833729,
 author = {Ramane, H. S. and Parvathalu, B. and Ashoka, K. and Pirzada, S.},
 title = {On families of graphs which are both adjacency equienergetic and distance equienergetic},
 fjournal = {Indian Journal of Pure \& Applied Mathematics},
 journal = {Indian J. Pure Appl. Math.},
 issn = {0019-5588},
 volume = {55},
 number = {1},
 pages = {198--209},
 year = {2024},
 language = {English},
 doi = {10.1007/s13226-022-00355-1},
 zbMATH = {7833729},
 Zbl = {1536.05305}
}

@article{zbMATH05664447,
 author = {Stevanovi{\'c}, Dragan and Indulal, Gopalapillai},
 title = {The distance spectrum and energy of the compositions of regular graphs},
 fjournal = {Applied Mathematics Letters},
 journal = {Appl. Math. Lett.},
 issn = {0893-9659},
 volume = {22},
 number = {7},
 pages = {1136--1140},
 year = {2009},
 language = {English},
 doi = {10.1016/j.aml.2008.11.007},
 zbMATH = {5664447},
 Zbl = {1179.05040}
}

@article{zbMATH07833946,
 author = {Tang, Liang Wen and Liu, Juan and Qin, Shuangwei},
 title = {Construction of equienergetic trees},
 fjournal = {MATCH - Communications in Mathematical and in Computer Chemistry},
 journal = {MATCH Commun. Math. Comput. Chem.},
 issn = {0340-6253},
 volume = {91},
 number = {2},
 pages = {485--488},
 year = {2024},
 language = {English},
 doi = {10.46793/match.91-2.485T},
 zbMATH = {7833946},
 Zbl = {1536.92184}
}

@article{zbMATH06221129,
 author = {Indulal, Gopalapillai},
 title = {Spectrum of two new joins of graphs and infinite families of integral graphs},
 fjournal = {Kragujevac Journal of Mathematics},
 journal = {Kragujevac J. Math.},
 issn = {1450-9628},
 volume = {36},
 number = {1},
 pages = {133--139},
 year = {2012},
 language = {English},
 zbMATH = {6221129},
 Zbl = {1289.05296}
}

@article{zbMATH08041341,
 author = {Shamsher, Tahir and Pirzada, S.},
 title = {Eigenvalues of the generalized subdivision graph with applications to graph energy},
 fjournal = {Discrete Mathematics, Algorithms and Applications},
 journal = {Discrete Math. Algorithms Appl.},
 issn = {1793-8309},
 volume = {16},
 number = {8},
 pages = {14},
 note = {Id/No 2350115},
 year = {2024},
 language = {English},
 doi = {10.1142/S179383092350115X},
 zbMATH = {8041341},
 Zbl = {1562.05233}
}

@article{zbMATH08090417,
 author = {Manjunatha, B. J. and Rakshith, B. R. and Veeresha, R. G.},
 title = {Distance equienergetic graphs of diameter 4},
 fjournal = {Examples and Counterexamples},
 journal = {Ex. Counterex.},
 issn = {2666-657X},
 volume = {7},
 pages = {3},
 note = {Id/No 100184},
 year = {2025},
 language = {English},
 doi = {10.1016/j.exco.2025.100184},
 zbMATH = {8090417},
 Zbl = {1570.05080}
}

@article{zbMATH07249081,
 author = {Ramane, H. S. and Patil, Daneshwari and Parvathalu, B. and Ashoka, K.},
 title = {Construction of distance equienergetic graphs through generalized composition},
 fjournal = {Journal of Advanced Mathematical Studies},
 journal = {J. Adv. Math. Stud.},
 issn = {2065-3506},
 volume = {13},
 number = {1},
 pages = {35--41},
 year = {2020},
 language = {English},
 zbMATH = {7249081},
 Zbl = {1451.05144}
}

@article{zbMATH06319411,
 author = {Aouchiche, Mustapha and Hansen, Pierre},
 title = {Distance spectra of graphs: a survey},
 fjournal = {Linear Algebra and its Applications},
 journal = {Linear Algebra Appl.},
 issn = {0024-3795},
 volume = {458},
 pages = {301--386},
 year = {2014},
 language = {English},
 doi = {10.1016/j.laa.2014.06.010},
 zbMATH = {6319411},
 Zbl = {1295.05093}
}

@article{zbMATH07090170,
 author = {Indulal, Gopalapilla and Scaria, Deena C. and Liu, Xiaogang},
 title = {The distance spectrum of the subdivision vertex join and subdivision edge join of two regular graphs},
 fjournal = {DML. Discrete Mathematics Letters},
 journal = {DML, Discrete Math. Lett.},
 issn = {2664-2557},
 volume = {1},
 pages = {36--41},
 year = {2019},
 language = {English},
 url = {www.dmlett.com/archive/DML19_v1_p3641.pdf},
 zbMATH = {7090170},
 Zbl = {1438.05089}
}

@article{brouwer2011spectra,
  title= "{Spectra of graphs}",
  author= "Brouwer, Andries E and Haemers, Willem H",
  year= "2011",
  journal= "Springer Science \& Business Media"
}
\end{document}